\documentclass[runningheads]{llncs}
\usepackage{graphicx}
\usepackage[utf8]{inputenc}
\usepackage[english]{babel}
\usepackage[utf8]{inputenc} 
\usepackage[bookmarks=false,hyperfootnotes=false]{hyperref}
\hypersetup{
			colorlinks=true,
			linkcolor=blue,
			anchorcolor=black,
			citecolor=green,
			urlcolor=blue
			}
\usepackage{tipa}
 \expandafter\let\csname equation*\endcsname\relax
  \expandafter\let\csname endequation*\endcsname\relax
\usepackage{amsmath,amsfonts,amssymb}

\usepackage{amsthm}
\usepackage{mathtools}

\newcommand{\RR}{\mathbb{R}}

\newcommand{\Kcal}{\mathcal{K}}

\newcommand{\Wcal}{\mathcal{W}}

\newcommand{\Ibold}{\mathbf{I}}

\newcommand{\Mbold}{\mathbf{M}}

\newcommand{\Sbold}{\mathbf{S}}

\newcommand{\abold}{{\bm{a}}}

\newcommand{\ebold}{{\bm{e}}}

\newcommand{\ubold}{{\bm{u}}}

\newcommand{\xbold}{{\bm{x}}}
\newcommand{\ybold}{{\bm{y}}}

\newcommand{\zbold}{{\bm{z}}}

\usepackage{bm} 

\newcommand{\lambdabold}{\bm{\lambda}}
\newcommand{\Psibold}{\bm{\Psi}}

\DeclareMathOperator*{\argmin}{argmin}

\AtBeginDocument{
  
  \let\argmin\relax

  \let\prox\relax

  \DeclareMathOperator*{\argmin}{arg\,min}

  \DeclareMathOperator{\prox}{\mathbf{prox}}

}

\usepackage{bm}
\newcommand{\grad}{{\bm{\nabla}}}

\newcommand{\blank}{\,{\cdot}\,}

\usepackage{tabularx}
\usepackage{booktabs}
\usepackage{multirow}
\usepackage{algorithmicx}
\usepackage[titlenumbered,ruled,vlined]{algorithm2e}
\usepackage{algpseudocode,float}
\algnewcommand\Input{\textbf{Input: }}
\algnewcommand\Package{\textbf{Package: }}
\algnewcommand\Parameters{\textbf{Parameters: }}
\algnewcommand\Output{\textbf{Output: }}

\usepackage{xparse}
   \DeclareDocumentCommand\TDVM{mmo}
 {%
  \IfValueTF{#3}
   {\mathrm{{TDV}}_{{#2}}^{#1}({#3})}
   {\mathrm{{TDV}}_{{#2}}^{#1}}%
 }

\DeclareDocumentCommand\LL{mo}
 {%
  \IfValueTF{#2}
   {\mathrm{{L}}^{#1}(#2)}
   {\mathrm{{L}}^{#1}}%
 }
 
\newcommand{\T}{\mathrm{T}} 

\newcommand{\diff}{\mathop{}\mathrm{d}}

\newcommand{\norm}[1]{\left\lVert #1 \right\rVert}

\usepackage{hyperref}

\begin{document}

\title{Total Directional Variation for Video Denoising
\thanks{
\tiny
SP acknowledges UK EPSRC grant EP/L016516/1 for the CCA DTC.
CBS acknowledges support from Leverhulme Trust project on Breaking the non-convexity barrier, EPSRC grant Nr.\ EP/M00483X/1, the EPSRC Centre EP/N014588/1, the RISE projects CHiPS and NoMADS, the CCIMI and the Alan Turing Institute.
}
}
%
%
\author{
Simone Parisotto\inst{1} \and Carola-Bibiane Sch\"onlieb\inst{2}
}
\authorrunning{S. Parisotto \and C.-B. Sch\"onlieb}
%
\institute{
CCA, University of Cambridge, Wilberforce Road, Cambridge, CB3 0WA, UK\\
\email{sp751@cam.ac.uk}
\and
DAMTP, University of Cambridge, Wilberforce Road, Cambridge, CB3 0WA, UK\\
\email{cbs31@cam.ac.uk} 
}

\maketitle              
\begin{abstract}
In this paper we propose a variational approach for video denoising, based on a total directional variation (TDV) regulariser proposed in \cite{ParMasSch18applied,ParMasSch18analysis} for image denoising and interpolation. In the TDV regulariser, the underlying image structure is encoded by means of weighted derivatives so as to enhance the anisotropic structures in images, e.g.\ stripes or curves with a dominant local directionality. 
For the extension of TDV to video denoising, the space-time structure is captured by the volumetric structure tensor guiding the smoothing process.
We discuss this and present our whole video denoising workflow. 
The numerical results are compared with some state-of-the-art video denoising methods.

\end{abstract}

\section{Introduction}
Video denoising refers to the task of removing noise in digital videos. 
Compared to image denoising, video denoising is usually a more challenging task due to the computational cost in processing large data and the redundancy of information, i.e.\ the expected similarity between two consecutive frames that should be inherited by the denoised video. 
A straightforward approach to video denoising is to denoise each frame of the video independently, by using the broad literature on image denoising methods, see e.g.\ 
\cite{bilateral,ROF,buadesimagedenoising,wavelet4,NL-means,BM3D,NL-Bayes,BreKunPoc2010,Papafitsoros2014,directionaltv,ParMasSch18applied}.
Computational cost is then stratified across image frames by sequentially processing them, which is seen as an advantage. However, a significant disadvantage of this frame-by-frame processing is the appearance of flickering artefacts and post-processing motion compensation step may be required \cite{motioncompensation,noisereductionreview}.

In recent years different approaches have been proposed for solving the video denoising problem: we refer to the introduction of \cite{AriasMorel2018} for an extensive survey. 
Notably, patch-based approaches are usually considered among the most promising video denoising methods in that they are able to achieve qualitatively good denoising results.
For example, V-BM3D is the 3D extension of the BM3D collaborative filters \cite{VBM3D}: without inspecting the motion time-consistency, V-BM3D independently filters 2D patches resulting similar in the 3D spatio-temporal neighbourhood domain.
As mentioned in \cite{AriasMorel2018}, while generally receiving good denoising results, the problem of flickering still occurs in V-BM3D.
For this reason, authors of V-BM3D developed an extension, called V-BM4D, where the patch-similarity is explored along space-temporal blocks defined by a motion vector, see \cite{VBM4D}.
Similarly, in \cite{Buades} the authors propose to group patches via an optical flow equation based on \cite{ZacPocBis07} and implemented in \cite{PerMeiFac13}.
In these approaches, while the incorporation of motion helps to provide consistency in time, denoising results also suffer from the lack of accuracy in the estimated motion.
A possible way to avoid the motion estimation is to consider 3D rectangular patches so as to inherently model the 3D structure and motion in the spatio-temporal video dimensions, based on the fact that rectangular 3D patches are less repeatable than motion-compensated patches.
However, such approach is not efficient for uniform motion or homogeneous spatial patterns, cf.\ the discussion on this topic in \cite{AriasMorel2018}.
Motivated by this reasoning the authors of \cite{AriasMorel2018} introduce a Bayesian patch-based video denoising approach with rectangular 3D patches modelled as independent and identically distributed samples from an unknown \emph{a priori} distribution: then each patch is denoised by minimising the expected mean square error.
Other approaches in video denoising are the straightforward extension of the Rudin-Osher-Fatemi (ROF) model \cite{ROF} to 3D data, by using a spatio-temporal total variation (referred in the next as ROF 2D+t), the joint video denoising with the computation of the flow \cite{BurDirSch18} and CNN approaches \cite{DavEhrFacMor18}.

\paragraph{Scope of the paper.}
In this paper we propose an extension of the recently introduced \emph{total directional variation} (TDV) regulariser \cite{ParMasSch18applied,ParMasSch18analysis} for video denoising, via the  
following variational regularisation model:
\begin{footnotesize}
\begin{equation}
u^\star \in \argmin_u \left(\TDVM{}{}[u,\Mbold] + \frac{\eta}{2}\norm{u-u^\diamond}_2^2\right),
\label{eq:model}
\end{equation}
\end{footnotesize}
where $u^\star$ is the denoised video, $\Mbold$ is a weighting field that encodes directional features in two spatial and one temporal dimension, $\eta>0$ is the regularisation parameter and $u^\diamond$ is a given noisy video. 
The model \eqref{eq:model} will be made more precise in the next sections where we mainly focus on its discrete and numerical aspects.
In order to accommodate for spatial-temporal data, we consider here a modification of the TDV regulariser given in \cite{ParMasSch18applied,ParMasSch18analysis} that derives directionality in the temporal dimension.
Differently from the patch-based approach, we compute for each voxel the vector field of the motion, to be encoded as a weight in the TDV regulariser. 
With this voxel-based approach we will reduce the flickering artefact 
which appears in patch-based approaches due to the patch selection, especially in regions of smooth motion.
Results are presented for a variety of videos corrupted with Gaussian white noise.

\paragraph{Organisation of the paper.} 
This paper is organised as follows:
in Section \ref{sec: tdv for video} we describe the estimation of the vector fields, the TDV regulariser and the variational model to be minimised; in Section \ref{sec: discrete} we describe the optimisation method for solving the TDV video denoising problem and comment on the selection of parameters; in Section \ref{sec: results} we show denoising results on a selection of videos corrupted with Gaussian noise of varying strength.

\section{Total directional variation for video denoising}
\label{sec: tdv for video}
Let $\overline{u}:\Omega \times [1,\dots,T]\to \RR_+^C$ be a clean video and $\Omega$ a spatial, rectangular domain indexed by $\xbold=(x,y)$, 
with number of $T$ frames and $C$ colours.
Let $u^\diamond$ be a corrupted version of $\overline{u}$ in each space-time voxel $(\xbold,t)\in\Omega\times[1,\dots,T]$ by i.i.d.\ Gaussian noise $n$ of zero mean and (possibly known) variance $\varsigma^2 >0$:
\begin{equation}\small
u^\diamond(\xbold,t) 
= 
\overline{u}(\xbold,t) + n(\xbold,t), \quad\forall (\xbold,t)\in\Omega\times[1,\dots,T].
\end{equation}
In what follows, we propose to compute a denoised video $u^\star\approx\overline{u}$ by solving
\begin{equation}\small
u^\star\in \argmin_u \left( \TDVM{}{}[u,\Mbold] + \frac{\eta}{2} \norm{u-u^\diamond}_2^2\right),
\label{eq: relaxation}
\end{equation}
where $\TDVM{}{}[u,\Mbold]$ is the proposed total direction regulariser w.r.t.\ a weighting field $\Mbold$, both specified in the next sections, and $\eta>0$ a regularisation parameter.


\subsection{The directional information}
In order to capture directional information of $u$ in \eqref{eq: relaxation}, we eigen-decompose the two-dimensional structure tensor \cite{weickert1998anisotropic} in each coordinate plane.

To do so, we first construct the 3D structure tensor: let $\rho\geq\sigma>0$ be two smoothing parameters,
$K_\sigma,K_\rho$ be the Gaussian kernels of standard deviation $\sigma$ and $\rho$, respectively, and 
let $u_\sigma = K_\sigma\ast u$. 
Then the 3D structure tensor 
reads as
\begin{equation}\small
\renewcommand{\arraystretch}{1}
\Sbold 
:= 
K_\rho \ast (\grad u_\sigma \otimes \grad u_\sigma) = \begin{pmatrix}
u_{\sigma,\rho}^{x,x} & u_{\sigma,\rho}^{x,y} & u_{\sigma,\rho}^{x,t}\\
u_{\sigma,\rho}^{y,x} & u_{\sigma,\rho}^{y,y} & u_{\sigma,\rho}^{y,t}\\
u_{\sigma,\rho}^{t,x} & u_{\sigma,\rho}^{t,y} & u_{\sigma,\rho}^{t,t}
\end{pmatrix},
\label{eq: 3d volumetric tensor}
\end{equation}
where $\grad u_\sigma \otimes \grad u_\sigma = \grad u_\sigma \grad u_\sigma^T$, $u_{\sigma,\rho}^{p,q}:=K_\rho\ast(\partial_p u_\sigma \otimes \partial_q u_\sigma)$ for each $p,q\in\{x,y,t\}$. 

\begin{figure}[!htb]
\centering
\includegraphics[width=0.25\textwidth]{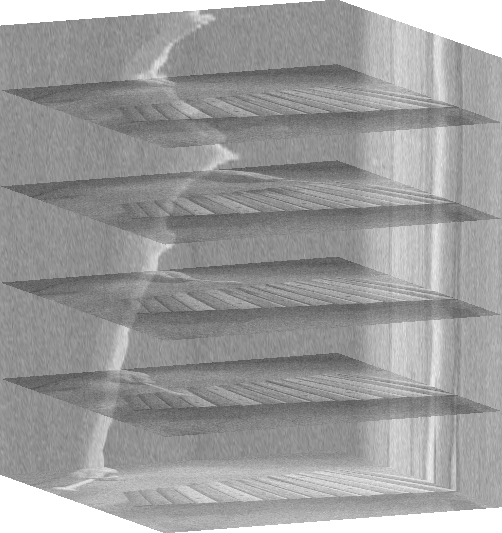}
\includegraphics[width=0.25\textwidth]{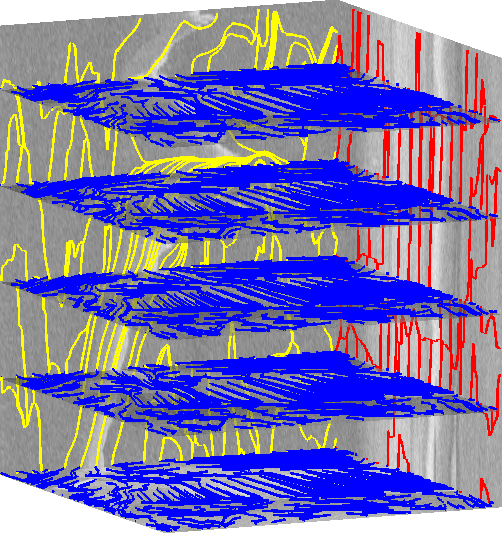}
\caption{
Left: grey-scale video \texttt{xylophone.mp4}, corrupted by Gaussian noise ($\varsigma=20$); right: streamlines of the weighting field with $\ebold_2$ (blue), $\ebold_4$ (red) and $\ebold_6$ (yellow).
}
\label{fig: streamlines}
\end{figure}

For a straightforward application to the TDV regulariser in \cite{ParMasSch18applied}, we extract the 2D sub-tensors of \eqref{eq: 3d volumetric tensor}, whose
eigen-decomposition encodes structural information in each of the coordinate frames spanned by $\{x,y\}$, $\{x,t\}$ and $\{y,t\}$:
\begin{equation}\renewcommand{\arraystretch}{1.1}\small
\begin{aligned}\small
\text{on coordinates $\{x,y\}$:}
\quad&
\Sbold^{x,y}=
\begin{pmatrix}
u_{\sigma,\rho}^{x,x} & u_{\sigma,\rho}^{x,y} \\
u_{\sigma,\rho}^{y,x} & u_{\sigma,\rho}^{y,y}
\end{pmatrix}
=
\lambda_1 (\ebold_{1}\otimes\ebold_{1})
+
\lambda_2 (\ebold_{2}\otimes\ebold_{2});
\\
\text{on coordinates $\{x,t\}$:}
\quad&
\Sbold^{x,t}=
\begin{pmatrix}
u_{\sigma,\rho}^{x,x} & u_{\sigma,\rho}^{x,t} \\
u_{\sigma,\rho}^{t,x} & u_{\sigma,\rho}^{t,t}
\end{pmatrix}
=
\lambda_3 (\ebold_{3}\otimes\ebold_{3})
+
\lambda_4 (\ebold_{4}\otimes\ebold_{4});
\\
\text{on coordinates $\{y,t\}$:}
\quad&
\Sbold^{y,t}=
\begin{pmatrix}
u_{\sigma,\rho}^{y,y} & u_{\sigma,\rho}^{y,t} \\
u_{\sigma,\rho}^{t,y} & u_{\sigma,\rho}^{t,t}
\end{pmatrix}
=
\lambda_5 (\ebold_{5}\otimes\ebold_{5})
+
\lambda_6 (\ebold_{6}\otimes\ebold_{6})
.
\end{aligned}
\label{eq: eigen-decompositions}
\end{equation}
For each $s\in\{1,\dots,6\}$, the eigenvector $\ebold_s=(e_{s,1},e_{s,2})$ has eigenvalue $\lambda_s$. 
The tangential directions in the 2D planes $\{x,y\},\{x,t\}$ and $\{y,t\}$ are $\ebold_2,\ebold_4,\ebold_6$, respectively, with
$\ebold_1,\ebold_3,\ebold_5$ the gradient directions, see Fig.\ \ref{fig: streamlines}.

From \eqref{eq: eigen-decompositions}, the ratios between the eigenvalues, called \emph{confidence}, measure the local anisotropy of the gradient on the slices within a certain neighbourhood:
\begin{footnotesize}
\begin{equation}
a^{x,y} = \dfrac{\lambda_2}{\lambda_1 + \varepsilon}, \quad
a^{x,t} = \dfrac{\lambda_4}{\lambda_3 + \varepsilon}, \quad
a^{y,t} = \dfrac{\lambda_6}{\lambda_5 + \varepsilon},\quad\text{with}\quad\varepsilon>0.
\label{eq: coherence}
\end{equation}
\end{footnotesize}
Here, $a^{x,y},a^{x,t},a^{y,t}\in[0,1]$ and
the closer to 0, the higher is the local anisotropy.
\vspace{-2em}

\vspace{-2.5em}
\subsection{The regulariser}
The TDV regulariser is composed of a gradient operator weighted by a tensor $\Mbold$, whose purpose is to smooth along selected directions. 
In view of the spatial-temporal data, we extend the natural gradient operator to the Cartesian planes $\{x,y\}$, $\{x,t\}$ and $\{y,t\}$. We will denote with $\widetilde{\grad}$ the concatenation of resulting 2-dimensional gradients. Further, 
we encode \eqref{eq: eigen-decompositions} and \eqref{eq: coherence} in $\Mbold$, leading to the \emph{weighted gradient} $\Mbold\widetilde{\grad}$ for the video function $u=u(x,y,t)$:
\begin{footnotesize}
\begin{align}
\Mbold \widetilde{\grad} \otimes u  
&=
\underbrace{
\begin{pmatrix}
a^{x,y} & 0 & 0 & 0 & 0 & 0 \\
0 & 1 & 0 & 0 & 0 & 0 \\
0 & 0 & a^{x,t}  & 0 & 0 & 0\\
0 & 0 & 0 & 1 & 0 & 0\\
0 & 0 & 0 & 0 & a^{y,t}  & 0\\
0 & 0 & 0 & 0 & 0 & 1
\end{pmatrix}
\begin{pmatrix}
e_{1,1} & e_{1,2} & 0 & 0 & 0 & 0 \\
e_{2,1} & e_{2,2}  & 0 & 0 & 0 & 0 \\
0 & 0 & e_{3,1} & e_{3,2} & 0 & 0\\
0 & 0 & e_{4,1} & e_{4,2} & 0 & 0\\
0 & 0 & 0 & 0 & e_{5,1} & e_{5,2}\\
0 & 0 & 0 & 0 & e_{6,1} & e_{6,2}
\end{pmatrix}
}_{\Mbold}
\underbrace{
\begin{pmatrix}
\partial_x\\
\partial_y\\
\partial_x\\
\partial_t\\
\partial_y\\
\partial_t
\end{pmatrix}
}_{\widetilde{\grad}}
\otimes u
\label{eq: first choice weights}\\
&=
\begin{pmatrix}
a^{x,y}\grad_{\ebold_1}^{x,y} u, &
\grad_{\ebold_2}^{x,y} u,&
a^{x,t}\grad_{\ebold_3}^{x,t} u,&
\grad_{\ebold_4}^{x,t} u,&
a^{y,t}\grad_{\ebold_5}^{y,t} u,&
\grad_{\ebold_6}^{y,t} u
\end{pmatrix}^\T.
\label{eq: explicit directional derivative}
\end{align}
\end{footnotesize}
Note that $\Mbold$ is computed once from the noisy input $u^\diamond$.
For a fixed frame $\{p,q\}$ with $p,q\in\{x,y,t\}$ and direction $\zbold=(z_1,z_2)$ the gradient $\grad_\zbold^{p,q} u=\partial_p u \cdot z_1 + \partial_q u\cdot z_2$ is the directional derivative of $u$ along $\zbold$ w.r.t.\ the frame $\{p,q\}$.
See \cite[Fig. 3.12]{ParisottoThesis} for more details about this choice.
With this notation in place, we consider the \emph{total directional variation} (TDV) regulariser,
\begin{footnotesize}
\begin{equation}
\TDVM{}{}[u,\Mbold]
= 
\sup_{\Psibold} \left\{
\int_\Omega 
 (\Mbold\widetilde{\grad}\otimes u) \cdot \Psibold \diff \xbold\,\Big\lvert\,\text{for all suitable test functions } \Psibold
 \right\}.
 \label{eq: modified regulariser}
\end{equation}
\end{footnotesize}
By plugging \eqref{eq: explicit directional derivative} into \eqref{eq: modified regulariser} we reinterpret \eqref{eq: modified regulariser} as a penalisation of the rate of change along $\ebold_2,\ebold_4,\ebold_6$, with coefficients $a^{x,y},a^{x,t},a^{y,t}$ as bias in the gradient estimation. Note, that while in \cite{ParMasSch18applied} the $\TDVM{}{}$ regulariser has been proposed for a general order of derivatives, we consider here only a $\TDVM{}{}$ regulariser of first differential order. 
\subsection{Connections to optical flow}
Let $(x,y,t)\in\Omega\times[1,\dots T]$ be a voxel and $u(x,y,t)$ its intensity in the grey-scale video sequence $u$. 
If $u(x,y,t)$ is moved by a small increment $(\delta_x,\delta_y,\delta_t)$ between two frames, then the \emph{brightness constancy} constraint reads
\begin{equation}\label{eq:brightnessconst}
u(x,y,t) = u(x+\delta_x, y+\delta_y, t+\delta_t).
\end{equation}
If $u$ is sufficiently smooth, then the \emph{optical flow} constraint is derived \cite{LucasKanade1981,HornSchunck1981} as a linearisation of \eqref{eq:brightnessconst} with respect to a velocity field $\zbold$:
\begin{equation}\small
\grad u(x,y,t)^T\cdot \zbold
= 
0,
\quad
\text{for all } (x,y,t) \in \Omega\times [1,T].
\label{eq: opticalflow}
\end{equation}
For a specific field $\zbold=(\widetilde{\zbold},1)$ with $\widetilde{\zbold} = (z_1(x,y),z_2(x,y))$, Equation \eqref{eq: opticalflow} is equivalent to
\begin{equation}\small
-\partial_t u 
= 
\partial_x u \cdot z_1 + \partial_y u \cdot z_2
=
\grad_{\widetilde{\zbold}}^{x,y} u\quad\text{for all } (x,y,t) \in \Omega\times [1,T].
\label{eq: opticalflow 2}
\end{equation}
We can now re-write \eqref{eq: opticalflow} by means of the following velocity vector fields:
\begin{small}
\begin{equation}
\begin{aligned}
a^{x,y}&(e_{1,1},e_{1,2},1), & a^{x,t}&(e_{3,1},1,e_{3,2}), & a^{y,t}&(1,e_{5,1},e_{5,2}), \\
& (e_{2,1},e_{2,2},1), &
& (e_{4,1},1,e_{4,2}),&
& (1,e_{6,1},e_{6,2}),
\label{eq: velocity fields}
\end{aligned}
\end{equation}
\end{small}
leading to
\begin{footnotesize}
\begin{equation}
\begin{pmatrix}
-a^{x,y}\partial_t u\\
-\partial_t u\\
-a^{x,t}\partial_y u\\
-\partial_y u\\
-a^{y,t}\partial_x u\\
-\partial_x u
\end{pmatrix}
=
\begin{pmatrix}
a^{x,y} \partial_x u \cdot e_{1,1} + a^{x,y} \partial_y u \cdot e_{1,2} \\
\partial_x u \cdot e_{2,1} + \partial_y u \cdot e_{2,2} \\
a^{x,t}\partial_x u \cdot e_{3,1} + a^{x,t}\partial_t u \cdot e_{3,2} \\
\partial_x u \cdot e_{4,1} + \partial_t u \cdot e_{4,2} \\
a^{y,t}\partial_y u \cdot e_{5,1} + a^{y,t}\partial_t u \cdot e_{5,2} \\
\partial_y u \cdot e_{6,1} + \partial_t u \cdot e_{6,2}
\end{pmatrix}
=
\begin{pmatrix}
a^{x,y}\grad^{x,y}_{\ebold_1} u \\
\grad^{x,y}_{\ebold_2} u \\
a^{x,t}\grad^{x,t}_{\ebold_3} u \\
\grad^{x,t}_{\ebold_4} u \\
a^{y,t}\grad^{y,t}_{\ebold_5} u \\
\grad^{y,t}_{\ebold_6} u 
\end{pmatrix}.
\label{eq: optical flow explicit}
\end{equation}
\end{footnotesize}
Here, the right-hand side of \eqref{eq: optical flow explicit} encodes the components that we aim to penalise in \eqref{eq: explicit directional derivative}. Thus, the penalisation of \eqref{eq: explicit directional derivative} is equivalent to the penalisation of the left-hand side of \eqref{eq: optical flow explicit}, assumed \eqref{eq: opticalflow} holds with velocity fields in \eqref{eq: velocity fields}.
Note that the weights $a^{x,y},a^{x,t}$ and $a^{y,t}$ 
add a contribution in the direction of the gradients $\ebold_1,\ebold_3,\ebold_5$, respectively.

\subsection{The minimisation problem}
We aim to find the denoised video $u^\star$ from the noisy input video $u^\diamond$ by solving the $\TDVM{}{}-\LL{2}$ minimisation problem
\eqref{eq: relaxation}.
For the numerical optimisation of \eqref{eq: relaxation} we use a primal-dual scheme \cite{chambolle_pock_2016}. For this, we rewrite \eqref{eq: relaxation} as a saddle point problem for the operator $\Kcal:=\Mbold\widetilde{\grad}$, whose adjoint will be denoted by $\Kcal^\ast$. 
In what follows, we denote by $u$ the primal variable, $y$ the dual variable, $f^\ast$ the Fenchel conjugate of $f$, by $g$ the fidelity term and by $\sigma,\tau>0$ the dedicated parameters of the primal-dual algorithm, see \cite{chambolle_pock_2016} for more details on the primal-dual schemes in image processing and \cite{ParMasSch18applied} for their application to variational problems with TDV regulariser. 
The resulting saddle-point problem reads
\begin{small}
\begin{equation}
u^\star 
\in 
\argmin_u \max_{y}
\Big(\langle \Kcal u, y \rangle - \underbrace{\delta_{\{\norm{\blank}_{2,\infty}\leq 1\}}(y)}_{f^\ast(y)} + \underbrace{\frac{\eta}{2}\norm{u-u^\diamond}_2^2}_{g(u^\diamond)}
\Big).
\label{eq: saddle point}
\end{equation}
\end{small}
In the primal-dual algorithm solving \eqref{eq: saddle point} we need the proximal operators:
\begin{small}
\begin{equation}
\begin{aligned}
\prox_{\sigma f^\ast}(y) = \frac{y}{\max\{1,\norm{y}_2\}},\quad
\prox_{\tau g}(u) = u + (\Ibold + \tau\eta)^{-1}\tau\eta (u^\diamond -u ),
\end{aligned}
\end{equation}
\end{small}
where $\Ibold$ is the identity matrix.
Note that $g$ is uniformly convex, with convexity parameter $\eta$, so the dual problem is smooth. An accelerated version of the primal-dual algorithm can be used in this case, e.g.\ \cite[Alg. 2]{ChaPoc2011}, starting with $\tau_0,\sigma_0>0$ where $\tau_0\sigma_0 L^2\leq 1$ and $L^2$ is the squared operator norm, $L^2:=\norm{\Kcal}^2\leq 24$ (which holds in connection with the discretisation in \eqref{eq: central finite difference} and stepsize $h=1$).

\section{The discrete model}
\label{sec: discrete}
In the discrete model, 
$\Omega$ is a rectangular grid of size $M\times N$ and 
a video $\ubold$ is a volumetric data of size $M\times N\times T\times C$
(height$\times$width$\times$frames$\times$colours).
Here, we consider grey-scale videos ($C=1$) along the axes $(i,j,k)\in\Omega\times[1,T]$, with $i=1,\dots,M$, $j=1,\dots,N$ and $k=1,\dots,T$.
An extension to coloured videos is straightforward by processing each colour channel separately.
Here, a fixed $(i,j,k)\in\Omega\times[1,T]$ identifies a voxel in the gridded video domain, i.e.\ a small cube of size $h$ in each axis direction. Then, $\ubold_{i,j,k}:=\ubold(i,j,k)$ is the intensity in the voxel $(i,j,k)$ in the grey-scale video sequence $\ubold$. 
The noisy input video is denoted by $\ubold^\diamond$ as well as the other discrete vectorial quantities, namely $\abold^{1,2},\abold^{1,3},\abold^{2,3}$ and $\lambdabold_s$ for $s=1,\dots,6$.

\subsection{Discretisation of derivative operators and vector fields}
We describe a finite difference scheme on the voxels 
by introducing the discrete gradient operator $\grad:\RR^{M\times N\times T}\to\RR^{M\times N\times T\times 3}$, with $\grad = (\partial_1,\partial_2,\partial_3)$ defined via the central finite differences on half step-size and Neumann conditions as
\begin{equation}\small
\begin{aligned}
(\ubold^1)_{i,j,k} := (\partial_1 \ubold)_{i+0.5,j,k} 
&=
\begin{dcases}
\frac{\ubold_{i+1,j,k}-\ubold_{i,j,k}}{h}, &\text{if } i=1,\dots,M-1,\\
0 &\text{if } i=M;
\end{dcases}
\\
(\ubold^2)_{i,j,k} := (\partial_2 \ubold)_{i,j+0.5,k} 
&=
\begin{dcases}
\frac{\ubold_{i,j+1,k}-\ubold_{i,j,k}}{h}, &\text{if } j=1,\dots,N-1,\\
0 &\text{if } j=N;
\end{dcases}
\\
(\ubold^3)_{i,j,k} := (\partial_3 \ubold)_{i,j,k+0.5} 
&=
\begin{dcases}
\frac{\ubold_{i,j,k+1}-\ubold_{i,j,k}}{h}, &\text{if } k=1,\dots,T-1,\\
0 &\text{if } k=T.
\end{dcases}
\end{aligned}
\label{eq: central finite difference}
\end{equation}
\begin{remark}
While $\ubold$ lies at the vertices of the discrete grid, $\grad \ubold$ lies on its edges. Thus, \eqref{eq: central finite difference} is advantageous for local anisotropy since it has sub-pixel precision and a more compact stencil radius than the classical forward scheme.
\end{remark}

In \eqref{eq: first choice weights}, $\widetilde{\grad}:\RR^{M\times N\times T}\to\RR^{M\times N\times T\times 6}$ acts on $\ubold$ as follows:
\begin{footnotesize}
\begin{equation}
\left(\widetilde{\grad} \otimes \ubold\right)_{i,j,k} :=
\begin{pmatrix}
\ubold^1, &
\ubold^2, &
\ubold^1, &
\ubold^3, &
\ubold^2, &
\ubold^3 
\end{pmatrix}_{i,j,k}^\T.
\end{equation}
\end{footnotesize}
Any field $\ebold_s$ with $s=1,\dots,6$ and confidence $a^{1,2}$, $a^{1,3}$, $a^{2,3}$ will be discretised in the cell centres $(i+0.5,j+0.5,k+0.5)$ of the discrete grid domain. The weighting multiplication in \eqref{eq: first choice weights} is performed via an intermediate averaging interpolation operator $\Wcal:\RR^{M\times N\times T \times 6}\to\RR^{(M-1)\times (N-1)\times (T-1) \times 6}$ that avoids artefacts due to the grid offset: this gives
$
\Mbold\Wcal\widetilde{\grad}:\RR^{M\times N\times T}\to \RR^{(M-1)\times(N-1)\times(T-1)\times 6}.
$

\subsection{TDV for video denoising}
The TDV-based workflow consists of two steps, with pseudo-code in Alg.\ \ref{alg: TDV for video denoising}. 
The first one computes the directions via the eigen-decomposition in \eqref{eq: eigen-decompositions} while the second one is the primal-dual algorithm \cite[Alg. 2]{ChaPoc2011}, whose stopping criterion is the root mean square difference between two consecutive dual variable iterates.

\begin{algorithm}[tb]\noindent 
\SetAlgoLined\footnotesize
\caption{TDV for video denoising}
\label{alg: TDV for video denoising}
\SetKwData{maxiter}{maxiter}
\SetKwData{tol}{tol}
\SetKwData{prox}{prox}
\SetKwData{colourchannel}{colour channel}
\SetKwInOut{Input}{Input}
\SetKwInOut{Output}{Output}
\SetKwInOut{Parameters}{Parameters}
\SetKwProg{Fn}{Function}{:}{}
\SetKwFunction{computeST}{compute\_3D\_structure\_tensor}
\SetKwFunction{computeD}{compute\_derivative\_operator}
\SetKwFunction{computeM}{compute\_weights}
\SetKwFunction{computeK}{compute\_K\_and\_adjoint}
\SetKwFunction{eigendecomposition}{eigendecomposition}
\SetKwFunction{computeanisotropy}{compute\_anisotropy}
\SetKwFunction{primaldual}{primal\_dual}
\SetKwFunction{TDVvideodenoising}{TDV\_video\_denoising}
\Input{A video $\ubold^\diamond\in[0,255]$ ($M\times N\times 1 \times T$), $\varsigma\in[0,255]$.}
\Output{the denoised video $\ubold$;}
\Parameters{for the primal-dual $\maxiter,\tol$; for the variational model: $(\sigma,\rho,\eta)$. 
}
\BlankLine
\Fn{\TDVvideodenoising}{
\BlankLine
\tcp{Compute operators for the weighted derivative}
$[\partial_1,\partial_2,\partial_3]$ = \computeD ($M,N,T$) \;
$\Sbold$ = \computeST($\ubold^\diamond,\sigma,\rho$)\;
$[\ebold_1,\ebold_2,\ebold_3,\ebold_4,\ebold_5,\ebold_6,\lambdabold_1,\lambdabold_2,\lambdabold_3,\lambdabold_4,\lambdabold_5,\lambdabold_6]$ = \eigendecomposition($\Sbold$)\;
$[\abold^{1,2},\abold^{1,3},\abold^{2,3}]$ = \computeanisotropy ($\lambdabold_1,\lambdabold_2,\lambdabold_3,\lambdabold_4,\lambdabold_5,\lambdabold_6$)\;
$\Mbold$ = \computeM ($\abold^{1,2},\abold^{1,3},\abold^{2,3},\ebold_1,\ebold_2,\ebold_3,\ebold_4,\ebold_5,\ebold_6$)\;
\BlankLine
\tcp{Proximal operators, adjoints and primal-dual from \cite{ChaPoc2011}}
$[\Kcal,\Kcal^\ast]$ = \computeK ($\Mbold,\partial_1,\partial_2,\partial_3$)\;
$\prox_{f^\ast}$ = @($\ybold$) $\ybold./\max\{1,\norm{\ybold}_2\}$ \;
$\prox_{g}$ = @($\ubold,\tau$) $\ubold + (\Ibold + \tau\eta)^{-1}\tau\eta (\ubold^\diamond - \ubold )$\;
$\ubold$ = \primaldual($\ubold^\diamond,\Kcal,\Kcal^\ast,\prox_{f^\ast},\prox_g,\maxiter,\tol$)\;
\BlankLine
}
\Return
\end{algorithm}

\section{Results}
\label{sec: results}
In this section we discuss the numerical results for video denoising obtained with Alg.\ \ref{alg: TDV for video denoising}.
Considered videos have been taken from a benchmark video dataset\footnote{\scriptsize 
Videos are freely available: 
\emph{Salesman} and \emph{Miss America} 
at \url{www.cs.tut.fi/~foi/GCF-BM3D}
\emph{Xylophone} in MATLAB; \emph{Water} (re-scaled, grey-scaled and clipped, Jay Miller, CC 3.0) at \url{www.videvo.net/video/water-drop/477}; \emph{Franke}'s function (a synthetic surface moving on fixed trajectories: the coloured one changes with the \texttt{parula} colormap).
}
\footnote{\scriptsize Results are available at \url{http://www.simoneparisotto.com/TDV4videodenoising}.}. 
Each video has values in $[0,255]$ corrupted with Gaussian noise.
We tested different noise levels with standard deviation $\varsigma=[10,20,35,50,70,90]$ without clipping the videos 
so as to conform to the observation model.

The quality of the denoised result $\ubold^\star$ is evaluated by the peak signal-to-noise ratio (PSNR) value w.r.t.\ a ground truth video $\overline{\ubold}$. 
The model requires the parameters $(\sigma,\rho,\eta)$ as input. Once provided, we solve the saddle-point minimisation problem in \eqref{eq: saddle point} 
via the accelerated primal-dual algorithm 
with $L^2=\norm{\Kcal}^2=24$, see \cite[Alg. 2]{ChaPoc2011}. Here the tolerance for the stopping criterion is fixed to $10^{-4}$ (on average reached in $300$ iterations). 
However, we experienced faster convergence and similar results with $L^2\ll 24$ and bigger tolerances, e.g.\ $10^{-3}$.

\subsection{Selection of parameters}
\label{sec: gradient ascend}
In the model, $\ubold^\star$ is sensitive to the choice of both $(\sigma,\rho)$ for the vector fields, and the regularisation parameter $\eta$ that is chosen according to the noise level. Choosing those parameters by a trial and error approach is computationally expensive and the best parameters may differ, even for videos with the same noise level. In particular, the parameters $(\sigma,\rho)$ depend on structure in the data, e.g.\ flat regions versus motions versus small details. Therefore, a strategy for tuning them is needed.

\begin{figure}[tbh]
\centering
\includegraphics[width=0.49\textwidth]{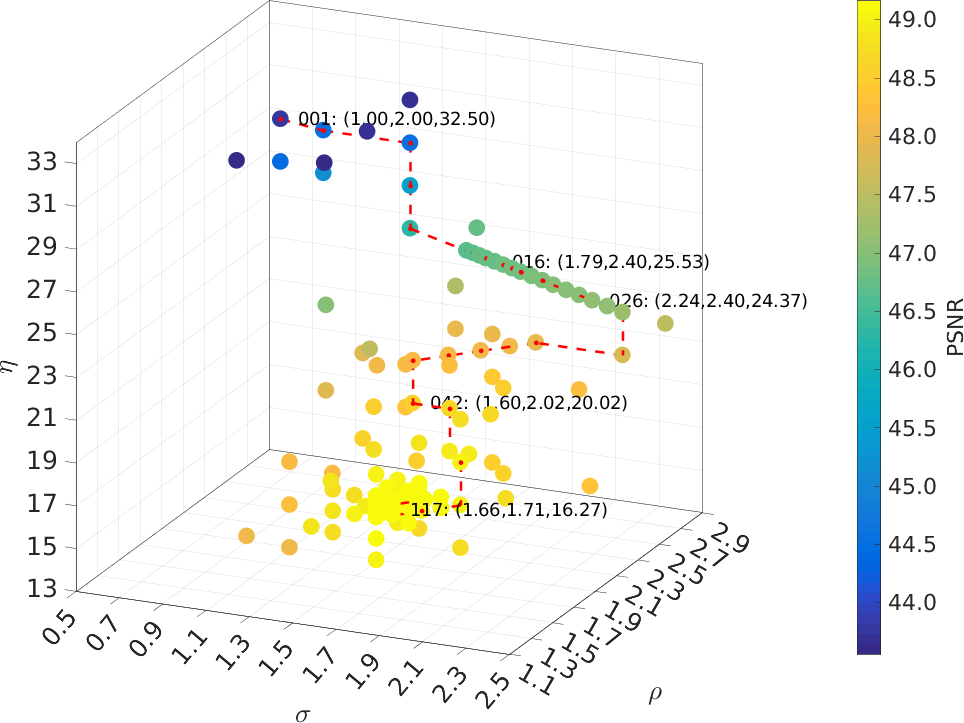}
\caption{
Line-search (\emph{Franke}, $\varsigma=10$): 
optimal trajectory (dashed red line);
PSNR values (coloured bullets). 
Optimal PSNR: 49.16, $117^\text{th}$ iteration, $(\sigma,\rho,\eta)=(1.66,1.71,16.27)$.}
\label{fig: linesearch function}
\end{figure}

To estimate appropriate values for $(\sigma,\rho,\eta)$ that render good results for a variety of videos we compute optimal parameters via line-search for maximising the PSNR for a small selection of video denoising examples for which the ground truth is available. The result of this optimisation is given in Table \ref{tab: PSNR franke}. 
For the line-search the parameters for the maximal PSNR values are computed iteratively, by applying Alg.\ \ref{alg: TDV for video denoising} for two different choices of $(\sigma,\rho,\eta)$ at a time, and subsequently adapt this parameter-set for the next iteration towards the ones in the neighbourhood of the one that returns a larger PSNR. In this search we constrain $\sigma\leq \rho$ \cite{weickert1998anisotropic}.
The line-search is stopped when, for the currently best parameters $(\sigma,\rho,\eta)$ all the other neighbours in a certain radius of distance report an inferior PSNR value. In Fig.\ \ref{fig: linesearch function} we show the trajectory of the parameters during this line-search for the \emph{Franke} video corrupted with Gaussian noise with $\varsigma=10$. We
observe that there exists a range of parameters in which the PSNR values are almost the same. 

By looking at the estimated parameters from the line-search approach in Table \ref{tab: PSNR franke}, we suggest the following rule of thumb for
their selection in Alg.\ \ref{alg: TDV for video denoising}:
\begin{equation}
\sigma = \rho = 3.2\eta^{-0.5} \quad\text{and}\quad\eta= 255\varsigma^{-1}.
\label{eq: quasi-optimal parameters}
\end{equation}

\subsection{Numerical results}
For the so-found optimal parameters we compare in Table \ref{tab: PSNR franke} 
the PSNR values achieved for our approach ($\TDVM{}{}$) with patch-based filters (V-BM3D v2.0 and V-BM4D v1.0, default parameters and \emph{normal-complexity profile}).
\begin{table}[tbh]\scriptsize
\centering
\caption{PSNR comparison (best in bold), with TDV parameters from line-search.}
    \label{tab: PSNR franke}
\begin{tabularx}{\textwidth}{Xc|c|c|c|c|c}
\toprule
Name ($M$, $N$, $C$, $T$) & $\varsigma$ & input & V-BM3D & V-BM4D & TDV ($\sigma,\rho,\eta$) & ROF 2D+t ($\eta$) \\
\midrule
                           & 10 & 28.13 & 45.99 & 46.90 & \textbf{49.16} (1.66,\,1.71,\,16.27) & 42.56 (16.27) \\
\emph{Franke grey-scale}    & 20 & 22.11 & 41.64 & 42.67 & \textbf{45.23} (2.00,\,2.00,\,08.10) & 38.18 (08.10)\\
$(120,120,1,120)$          & 35 & 17.25 & 38.63 & 39.34 & \textbf{41.89} (2.40,\,2.40,\,04.70) & 34.59 (04.70) \\
                           & 50 & 14.15 & 36.37 & 37.17 & \textbf{39.64} (2.70,\,2.70,\,03.30) & 32.30 (03.30) \\
                           & 70 & 11.23 & 30.60 & 35.03 & \textbf{37.44} (3.00,\,3.00,\,02.45) & 30.45 (02.45)\\
\midrule
                           & 10 & 28.13 & 47.13 & 48.21 & \textbf{50.51} (1.89,\,1.92,\,16.59) & 44.10 (16.59) \\
\emph{Franke coloured}     & 20 & 22.11 & 42.96 & 43.97 & \textbf{46.46} (2.35,\,2.35,\,08.35) & 39.93 (08.35) \\
$(120,120,3,120)$          & 35 & 17.25 & 40.18 & 40.47 & \textbf{42.97} (2.79,\,2.83,\,04.74) & 36.36 (04.74)\\
                           & 50 & 14.15 & 38.11 & 38.15 & \textbf{40.74} (3.13,\,3.17,\,03.45) &  34.41 (03.45) \\
                           & 70 & 11.23 & 31.72 & 35.90 & \textbf{38.62} (3.50,\,3.50,\,02.45) & 32.29 (02.45) \\
                           \midrule
                           & 10 & 28.13 & \textbf{37.30} & 37.12 & 35.24 (0.55,\,0.68,\,29.25) & 31.48 (29.25) \\
\emph{Salesman}            & 20 & 22.11 & \textbf{34.13} & 33.33 & 31.96 (0.70,\,0.75,\,13.93) &  28.16 (13.93) \\
$(288,352,1,050)$          & 35 & 17.25 & \textbf{30.79} & 30.20 & 29.36 (0.89,\,0.89,\,07.95) & 26.01 (07.95)\\
                           & 50 & 14.15 & 28.32 & \textbf{28.33} & 27.78 (1.05,\,1.06,\,05.45) & 24.78 (05.45) \\
                           & 70 & 11.23 & 24.55 & \textbf{26.68} & 26.34 (1.27,\,1.32,\,03.96) & 23.87 (03.96) \\
\midrule
                           & 10 & 28.13 & 43.83 & \textbf{44.68} & 43.13 (0.93,\,1.15,\,25.75) & 39.18 (25.75)\\
\emph{Water}               & 20 & 22.11 & 40.59 & \textbf{41.02} & 39.84 (1.18,\,1.35,\,12.60) & 35.94 (12.60)\\
$(180,320,1,120)$          & 35 & 17.25 & 37.75 & \textbf{37.90} & 37.14 (1.40,\,1.40,\,06.95) & 33.36 (06.95)\\
                           & 50 & 14.15 & 35.58 & \textbf{35.85} & 35.41 (1.61,\,1.65,\,04.80) & 31.83 (04.80)\\
                           & 70 & 11.23 & 30.11 & \textbf{33.87} & 33.78 (1.80,\,1.85,\,03.45) &  30.51 (03.45)\\
\bottomrule
\end{tabularx}
\end{table}

\begin{figure}[tbh]
\centering
\begin{minipage}[b]{0.45\textwidth}
\centering
\includegraphics[width=1\textwidth]{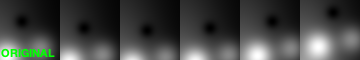}
\includegraphics[width=1\textwidth]{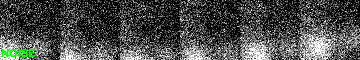}
\includegraphics[width=1\textwidth]{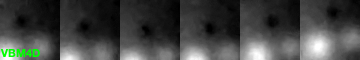}
\includegraphics[width=1\textwidth]{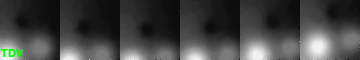}
\includegraphics[width=1\textwidth]{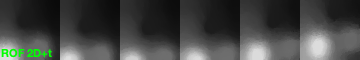}
\end{minipage}
\hfill
\begin{minipage}[b]{0.3\textwidth}
\centering
\includegraphics[width=1\textwidth]{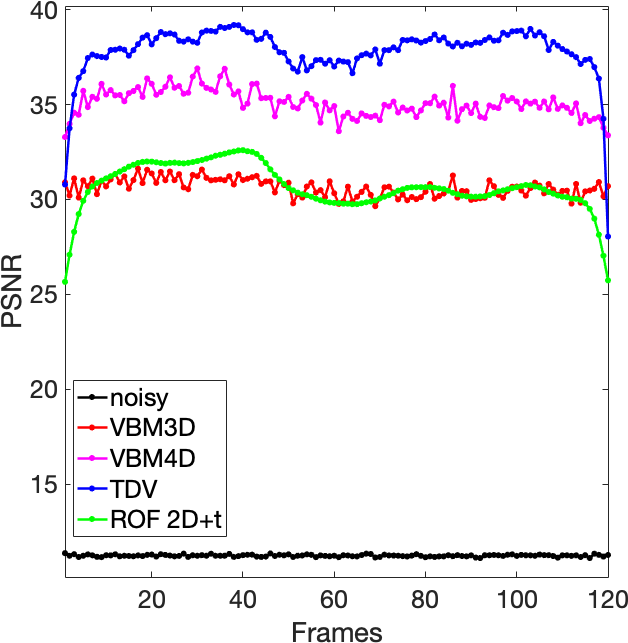}
\end{minipage}
\caption{
\emph{Franke}: frames $[10,20,30,40,50,60]$ for $\varsigma=70$ and PSNR comparison.
}
\label{fig: function}
\vspace{0.1em}
\begin{minipage}[b]{0.55\textwidth}
\centering
\includegraphics[width=1\textwidth]{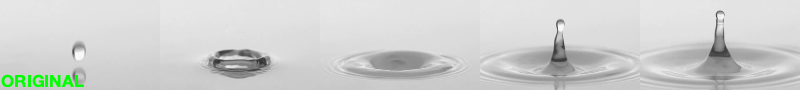}
\includegraphics[width=1\textwidth]{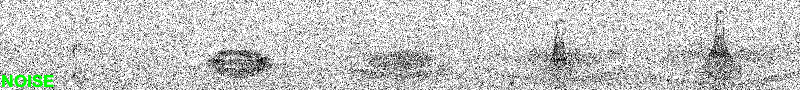}
\includegraphics[width=1\textwidth]{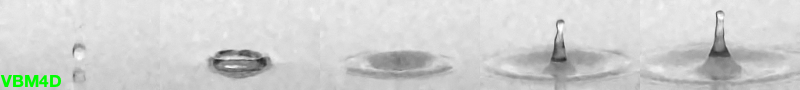}
\includegraphics[width=1\textwidth]{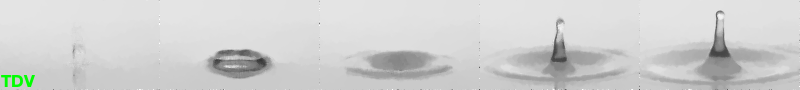}
\includegraphics[width=1\textwidth]{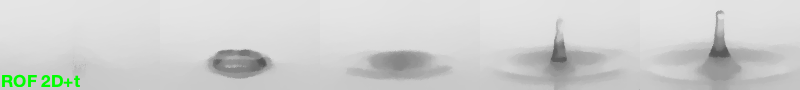}
\end{minipage}
\hfill
\begin{minipage}[b]{0.3\textwidth}
\centering
\includegraphics[width=1\textwidth]{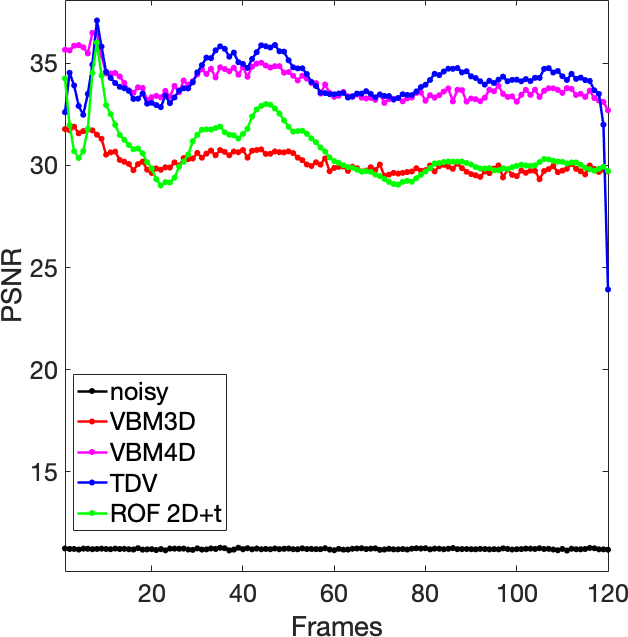}
\end{minipage}
\caption{
\emph{Water}: frames $[5,25,45,65,85]$ for $\varsigma=70$ and PSNR comparison.
}
\label{fig: water}
\end{figure}

In Figs.\ \ref{fig: function} and \ref{fig: water} the visual comparison is shown for selected frames of the \emph{Franke} and \emph{Water} videos (corrupted by noise with $\varsigma=70$). The time-consistency achieved by our approach is apparent in the frame-by-frame PSNR comparison.

\begin{figure}[tb]
\centering
\begin{minipage}[b]{0.55\textwidth}\centering
\includegraphics[width=1\textwidth]{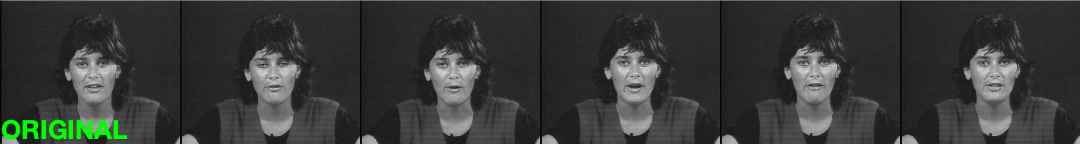}
\includegraphics[width=1\textwidth]{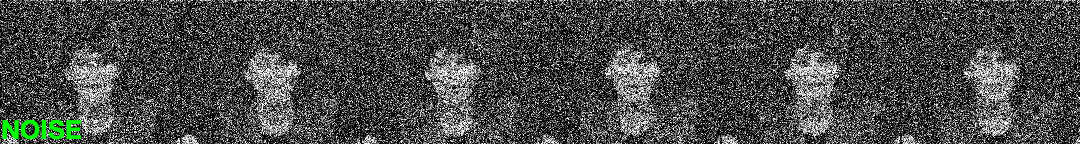}
\includegraphics[width=1\textwidth]{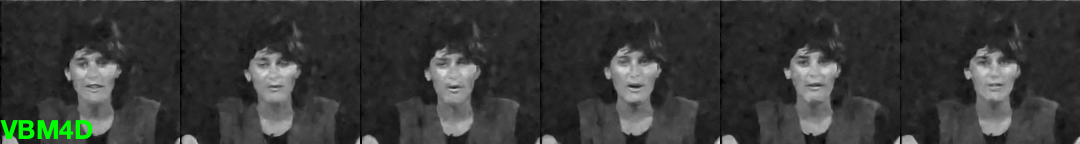}
\includegraphics[width=1\textwidth]{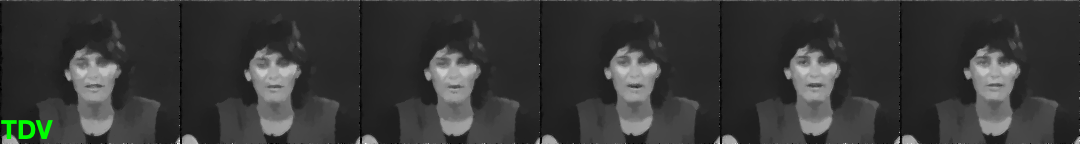}
\includegraphics[width=1\textwidth]{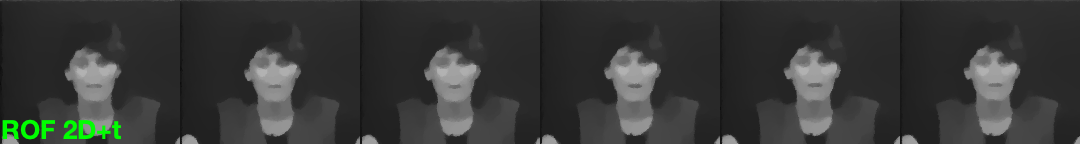}
\end{minipage}
\hfill
\begin{minipage}[b]{0.3\textwidth}\centering
\includegraphics[width=1\textwidth]{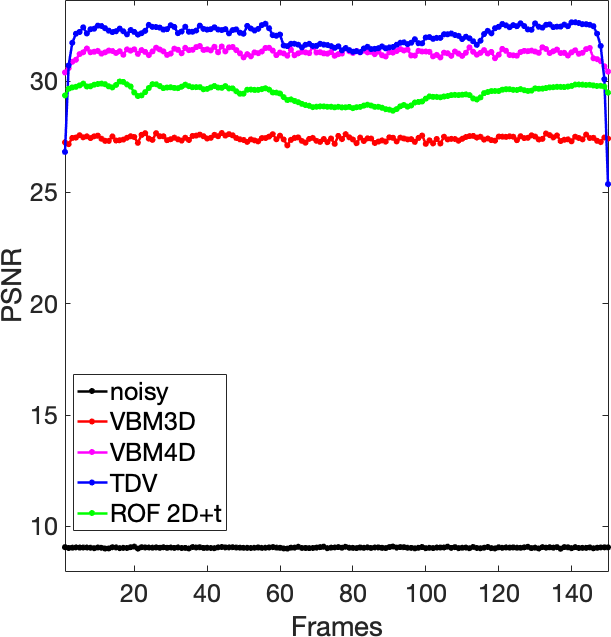}
\end{minipage}
\caption{
\emph{Miss America}: frames $[5,10,15,20,25,30]$ for $\varsigma=90$ and PSNR comparison.
}
\label{fig: miss}
\begin{minipage}[b]{0.55\textwidth}\centering
\includegraphics[width=1\textwidth]{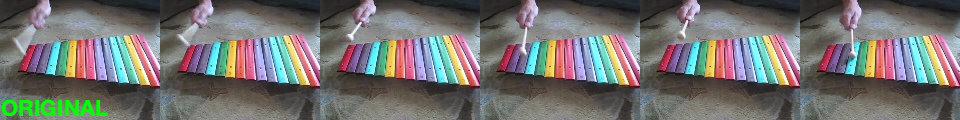}\\
\includegraphics[width=1\textwidth]{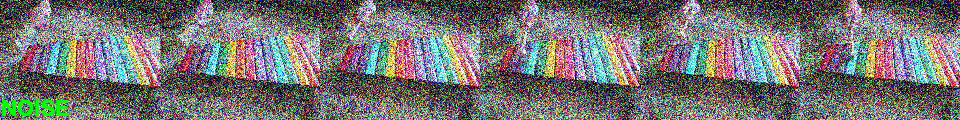}\\
\includegraphics[width=1\textwidth]{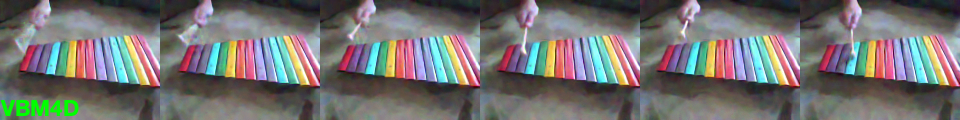}\\
\includegraphics[width=1\textwidth]{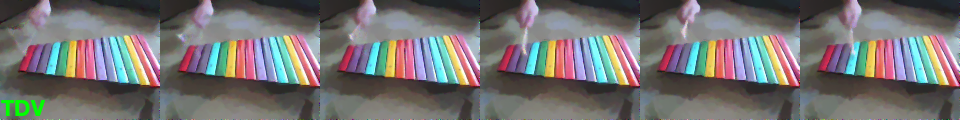}\\
\includegraphics[width=1\textwidth]{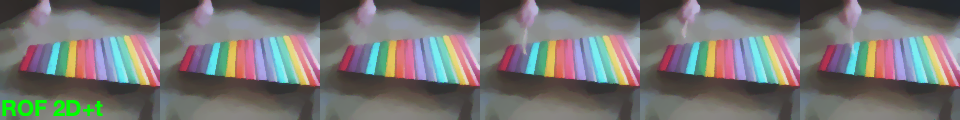}
\end{minipage}
\hfill
\begin{minipage}[b]{0.3\textwidth}\centering
\includegraphics[width=1\textwidth]{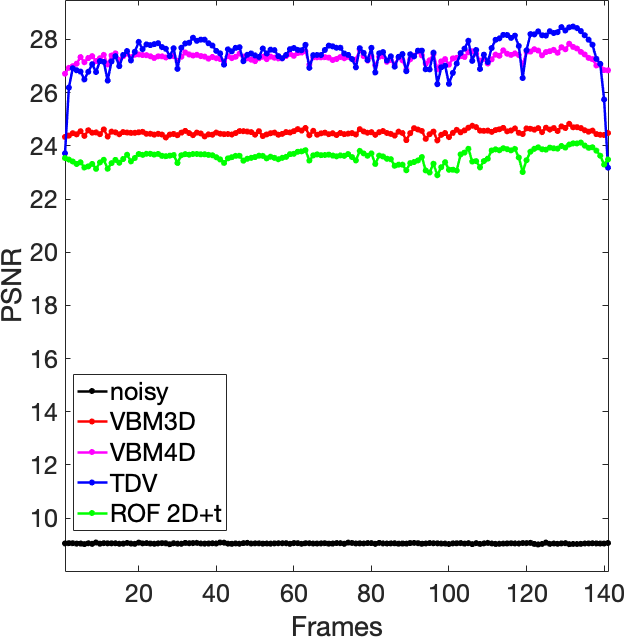}
\end{minipage}
\caption{
\emph{Xylophone}: frames $[5,10,15,20,25,30]$ for $\varsigma=90$ and PSNR comparison.
}
\label{fig: xylophone}
\end{figure}

Video denoising results that use the quasi-optimal parameters computed with \eqref{eq: quasi-optimal parameters} 
are reported in Table \ref{tab: PSNR quasi-optimal}: selected frames of videos corrupted with a high noise level of $\varsigma=90$ are shown in Figs.\ \ref{fig: miss}
and \ref{fig: xylophone}, with frame-by-frame PSNR values. 

\begin{table}[tbh]\scriptsize
\caption{PSNR comparison (best in bold), with quasi-optimal TDV parameters.}
\label{tab: PSNR quasi-optimal}
\begin{tabularx}{\textwidth}{Xc|c|c|c|c|c}
\toprule
Name ($M$, $N$, $C$, $T$) & $\varsigma$ & input & V-BM3D & V-BM4D & TDV ($\sigma,\rho,\eta$) & ROF 2D+t ($\eta$)\\
\midrule
                           & 10 & 28.13 & \textbf{39.64} & 39.93 & 39.25 (0.63,\,0.63,\,25.50) & 36.93 (25.50)\\
\emph{Miss America}        & 20 & 22.11 & \textbf{37.95} & 37.78 & 37.28 (0.90,\,0.90,\,12.75) & 34.60 (12.75) \\
$(288,360,1,150)$          & 35 & 17.25 & \textbf{36.03} & 35.77 & 35.44 (1.19,\,1.19,\,07.29) & 32.72 (07.29)\\
                           & 50 & 14.15 & \textbf{34.19} & 34.26 & 34.14 (1.42,\,1.42,\,05.10) & 31.47 (05.10)\\
                           & 70 & 11.23 & 28.86 & 32.64 & \textbf{32.85} (1.68,\,1.68,\,03.64) & 30.29 (03.64) \\
                           & 90 & 09.05 & 27.42 & 31.27 & \textbf{31.87} (1.90,\,1.90,\,02.83) & 29.42 (02.83)\\
\midrule
                           & 10 & 28.13 & \textbf{37.82} & 37.49 & 35.96 (0.63,\,0.63,\,25.50) & 32.70 (25.50)\\
\emph{Xylophone coloured}  & 20 & 22.11 & \textbf{34.70} & 34.13 & 33.06 (0.90,\,0.90,\,12.75) & 29.57 (12.75)\\
$(240,320,3,141)$          & 35 & 17.25 & \textbf{32.06} & 31.65 & 30.93 (1.19,\,1.19,\,07.29) & 27.16 (07.29)\\
                           & 50 & 14.15 & 29.98 & \textbf{30.07} & 29.58 (1.42,\,1.42,\,05.10) & 25.72 (05.10)\\
                           & 70 & 11.23 & 25.89 & \textbf{28.51} & 28.32 (1.68,\,1.68,\,03.64) & 24.43 (03.64) \\
                           & 90 & 09.05 & 24.50 & 27.32 & \textbf{27.37} (1.90,\,1.90,\,02.83) & 23.57 (02.83) \\
\bottomrule
\end{tabularx}
\end{table}

\subsection{Discussion of results}
We compared our variational $\TDVM{}{}$ denoising approach with patch-based (V-BM3D/V-BM4D) and variational (ROF 2D+t) methods.
Patch-based methods are usually computationally faster than the variational approaches (including ours) but they tend to suffer from flickering and staircasing artefacts due to their patch-based nature. 
We experienced that our MATLAB code (not optimised for speed) is approximately $7\times$ slower than V-BM4D (C++ code with MEX interface) with \emph{normal-complexity profile}. 
Both quantitative (via PSNR) and qualitative results (visual inspection) are relevant indicators for video denoising.

From the PSNR values in Tables \ref{tab: PSNR franke} and \ref{tab: PSNR quasi-optimal} the TDV approach is comparable with the patch-based ones, with many single frames achieving higher PSNR value than the patch-based methods did.
Also, by changing the noise level,
the PSNR values are deteriorating less than with the patch-based methods, demonstrating the consistency of our approach.

Visual results confirm that the $\TDVM{}{}$ approach improves upon patch-based methods producing less flickering and stair-casing artefacts, especially when the motion is smooth due to the coherence imposed also along the time dimension.


\section{Conclusions}
In this paper, we proposed a variational approach with the total directional variation (TDV) regulariser for video denoising.
We extended the range of applications of $\TDVM{}{}$ regularisation from image processing as demonstrated in \cite{ParMasSch18applied} to videos.
We compared $\TDVM{}{}$ with some state of the art patch-based algorithms for video denoising and obtained comparable results especially for high level noises while reducing artefacts in regions with smooth large motion, where the patch-based approach shows some weakness.
We expect to improve further the results by refining the estimation of the anisotropic fields \cite{Buades} and by using higher-order derivatives in the TDV definition  \cite{ParMasSch18applied}. This is left for future research.

\bibliographystyle{splncs04}
\bibliography{biblio}

\end{document}